\documentclass{article} 

\usepackage{amsmath}
\usepackage{amssymb}

\usepackage{amsthm}
\usepackage{amscd}

\newenvironment{pf}{\begin{proof}}{\end{proof}}
\newtheorem{thm}{Theorem}[section]
\newtheorem{lem}[thm]{Lemma}
\newtheorem{cor}[thm]{Corollary}
\newtheorem{prop}[thm]{Proposition}

\newtheorem{defn}{Definition}

\def\C{{\mathbb C}}

\def\P{{\mathbb P}}
\def\Q{{\mathbb Q}}

\def\Z{{\mathbb Z}}

\def\ccI{{\cal I}}
\def\ccJ{{\cal J}}
\def\ccK{{\cal K}}

\def\ccO{{\cal O}}
\def\ccQ{{\cal Q}}

\def\ccS{{\cal S}}

\def\ccU{{\cal U}}
\def\ccV{{\cal V}}

\def\ccX{{\cal X}}

\def\da{\downarrow}

\def\iso{\cong}

\def\ilim{\varprojlim}

\def\lra{\longrightarrow}
\def\m{{\mathfrak m}}

\def\ra{\rightarrow}

\def\tensor{\otimes}

\def\operatorname#1{\mathop{\rm #1}\nolimits}

\def\Chow{\operatorname{Chow}}

\def\H{{\rm H}}

\def\Hom{\operatorname{Hom}}

\def\Pic{\operatorname{Pic}}

\def\Spec{\operatorname{Spec}}

\def\im{\operatorname{im}}

\newcommand{\mr}[1]{\rightarrow}

\def\f{{\varphi}}

\def\remark#1{}

\title{Small Contractions of Symplectic $4$-folds}

\author{Jan Wierzba and Jaros\l{}aw A. Wi\'sniewski }

\begin{document}

\maketitle
\noindent
{\bf Abstract:} \ \ \ We classify small contractions of
(holomorphically) symplectic 4-folds. 

\tableofcontents

\section{Introduction}

\begin{defn} A {\em symplectic manifold} is defined to be 
a complex algebraic or analytic manifold $X$ of dimension $2n$, 
which carries a (holomorphic) symplectic $2$-form 
$\sigma\in \H^0(X,\Omega_X^2)$, that is $\sigma$ is closed and
anywhere nondegenerate, which means $d\sigma=0$ and
$\sigma^{\wedge n}$ vanishes nowhere. \end{defn}

We note that the form $\sigma^{\wedge n}$ trivializes the
canonical bundle $K_X$ so if $X$ is compact and simply-connected then it
is a Calabi-Yau variety. The symplectic form $\sigma$ defines an
isomorphism $\hat\sigma: T_X\ra\Omega^1_X$ given by the formula
$\hat\sigma(v)(w)=\sigma(v,w)$; we will call it a $\sigma$-duality. A
reduced but possibly reducible subvariety $M\subset X$ is called
Lagrangian if any component $M$ has dimension $n$ and at any
smooth point $x$ of $M$ the form $\sigma$ is trivial on the tangent
space $T_xM$.

We want to understand the local structure of birational morphisms of
projective symplectic manifolds. For this purpose we introduce
the following definition.

\begin{defn} Let $\f:X\ra Y$ be a birational projective morphism (of
complex algebraic varieties or complex analytic spaces) where $X$ is a
symplectic manifold and where $Y$ is normal. We say that $\f:X\ra Y$ is a
{\em symplectic contraction}. \end{defn}

We note that $K_X=\f^*K_Y$ so $\f$ is a crepant (or log terminal)
contraction in the sense of the Minimal Model Program. In terms of the
Program crepant contractions form a natural extension of Fano-Mori
contractions for which $-K_X$ is $\f$-ample (for a discussion on these
see e.g.~\cite{AW-SC}).  As it will turn out, symplectic contractions
have some features which make them somewhat similar to Fano-Mori
contractions.

In what follows we will frequently shrink both the domain as well as the
target of the contraction in order to understand its local structure.
Sometimes we will stress it by referring to them as to {\sl local}
contractions. Here is an example of a symplectic contraction.
\medskip

\noindent {\bf Example.} Let $X=T^*\P^2=\Spec_{\P^2}(\bigoplus_{m\geq 0}
S^m T_{\P^2})$ be the cotangent bundle of $\P^2$. Then $X$ is a
smooth variety of dimension $4$ which carries a natural symplectic form.
The bundle $T_{\P^2}$ is very ample, and if we set  $Y=\Spec
\bigoplus_{m\geq 0} H^0(\P^2,S^m T_{\P^2})$ then $Y$ is  a variety
of dimension $4$ and we have a birational morphism $\f:X\ra Y$, defined
by the evaluation of global sections $$\bigoplus_{m\geq 0} H^0(\P^2,S^m
T_{\P^2})\ra \bigoplus_{m\geq 0} S^m T_{\P^2}$$ which contracts
precisely the zero section of $T^* \P^2$ to  a single point. The
morphism $\f:X\ra Y$ is a symplectic contraction and we  refer to it as
the {\sl collapsing of the zero section in the cotangent bundle of}
$\P^2$. 

\medskip

Let us note that conversely, if a symplectic 4-fold $X$ contains a
subvariety $E\cong \P^2$ then a formal neighborhood of $E$ is
isomorphic to the formal neighborhood of the zero section of $T^*\P^2$.
Indeed, it is easy to verify that $E$ is a Lagrangian subvariety with
normal bundle isomorphic to $T^*_{\P^2}=\Omega_{\P^2}$. Since for
$i>0$ we have the vanishing  $H^1(\P^2,T_{\P^2}\otimes
S^i(T_{\P^2}))=H^1(\P^2, \Omega_{\P^2}\otimes S^i(T_{\P^2}))=0$
then by Grauert and Hironaka-Ross (see \cite[3.33]{Mori82}, 
\cite[3.7]{AW-SC}) 
an analytic neighborhood is uniquely defined.

The variety $Y$ in our example is equal to the affine cone over the
projective variety $\P(T_{\P^2})\subset \P^7$ which is a hyperplane
section of the Segre embedding $\P^2\times\P^2\subset\P^8$. The morphism
$\f$ is a resolution of the vertex singularity which is the graph of a
rational morphism obtained from the projection $\P(T_{\P^2})
\subset\P^2\times\P^2\ra \P^2$. The two projections lead to two
non-equivalent resolutions $\f: X\ra Y$ and $\f': X'\ra Y$ which are
dominated by a simple blow-up of the vertex.

More generally, let $X$ be a symplectic manifold of dimension $4$ and
suppose $X$ contains a $\P^2$. Let $\beta:Z\ra X$ be the blow up of
$X$ along $\P^2$.  As we noted above, in an analytic neighborhood of
the exceptional divisor of $\beta$ the manifold $Z$ is isomorphic to
the cotangent bundle of $\P^2$ blown-up along zero section. Thus there
exists another blow-down map $\beta':Z\ra X'$, where $X'$ is again a
symplectic manifold with a symplectic form $\sigma'$ which coincides
with $\sigma$ outside of the exceptional locus of the transformation.
We note that although the above arguments are performed in analytic
category their algebraic counterpart holds whenever one assumes that
the $\P^2$ in question is an isolated positive-dimensional fiber of a
symplectic contraction $\f: X\ra Y$. Indeed, in such a case the
blow-down $\beta':Z\ra X'$ exists by the relative cone and contraction
theorems, see e.g. \cite[3.25]{KM} or \cite[4-2-1]{KMM}.

\begin{defn} The birational map $\phi:X - - \ra X'$ constructed above is
called the {\em Mukai flop}.  \end{defn}

The main result of the paper is the following.

\begin{thm}\label{big-thm} Let $\f:X\ra Y$ be a symplectic contraction
with $Y$ quasiprojective.
Suppose that $\dim X=4$ and that $\f$ is small (i.e. it does not contract
a divisor). Then $\f$ is local analytically isomorphic to the collapsing
of the zero section in the cotangent bundle of $\P^2$ and therefore it
admits a Mukai flop. \end{thm}

We note that in view of the preceding discussion the actual contents of
the theorem is as follows: if $\f$ is a small contraction of a
symplectic 4-fold then the exceptional locus of $\f$ is isomorphic to
$\P^2$. The assumption of $Y$ quasiprojective is probably not needed,
that is we believe that the target of $\f$ may be just an analytic space,
we put it here in order to use deformation results: \cite[Thm.~2.1]{Wierzba}
and \cite[Thm.~1.3]{Kaledin00}, 
which we need for the crucial observation \ref{hom-4-connected}.

\medskip

Burns, Hu and Luo observed in \cite{BHL} that the above theorem allows
to understand birational transformations of symplectic 4-folds.  
A complete proof of this result
will be given in a forthcoming paper by the first named author. 
Using 
\ref{big-thm} and the Minimal Model Program, 
in particular \cite[Prop. 2.7]{Corti}, and
standard arguments providing termination of log-flips, 
one gets the following.

\begin{thm}\label{thm-mukai} Let $\phi:X - - \ra X'$ be a birational
map of two smooth projective 4-folds which is isomorphism in codimension
1.  Suppose that $X$ is symplectic. Then $X'$ is symplectic as well
and $\phi$ can be factorised in a finite sequence of Mukai
flops. \end{thm}

The paper is organised as  follows. Firstly we recall results which we
use in  the course of the  proof of \ref{big-thm} and  restate them in
the  suitable context.  This concerns  results about  rational curves,
vanishings,  normal  surfaces  with  many rational  curves  and  local
existence of a vector field on a symplectic manifold. The actual proof
of  \ref{big-thm} is  divided into  four  steps: first  we reduce  the
arguments to the case when  $\f$ is an elementary contraction. Then we
use deformation of rational curves, as developed by Mori and Koll\'ar,
in order  to understand  1-dimensional singular locus  of $E$:  as the
result we prove that $E$  is irreducible and homeomorphic to $\P^2$ in
codimension 1.  Using this we are  able to produce locally  near $E$ a
line  bundle  $\ccO(1)$ for  which  we  apply Kawamata  base-point-free
arguments.   Finally  we  analyse  a  non-normal  del  Pezzo  surface,
appearing  in  a classification  list  by Reid,  which  is  left as  a
possible exception by the previous argument.

Recently (September 2001) we were informed that Cho, Miyaoka and
Shepherd-Barron have obtained a proof of a version of \ref{big-thm}
valid in arbitrary dimension, see \cite{CMSB}. Their method however is
substatially different from ours.

\bigskip

\noindent{\bf Thanks and acknowledgements.} We would like to thank
people by communicating with whom we benefited during the preparation
of this paper. The first named author thanks Nick Shepherd-Barron, who
was his scientific advisor in Cambridge, as well as Daniel Huybrechts,
Manfred Lehn and Dmitri Kaledin, whom he met in K\"oln. The second
named author would like to thank Stavros Papadakis, Luis Sola Conde
and Massimiliano Mella.

Our collaboration was made possible due to EC grant
HPRN-CT-2000-00099: the first named author was a post-doc fellow of
EAGER program in Warsaw.  The second named author would like to
acknowledge partial support of the Polish KBN.


\section{Preliminaries}

Let $\f:X\ra Y$ be a morphism of analytic spaces, which is a
symplectic contraction. The morphism $\f$ is then a crepant or
log-extremal contraction in the sense of the Minimal Model Program and
as such it admits several special properties which we want to discuss
in this preliminary section. For the fundamental results of the
Program we refer to \cite{KMM}, \cite{KM} and \cite{Kawamata88} in the
analytic set-up. When dealing with rational curves we will use the
language and notation consistent with \cite{Kollar}.

Let $E=\bigcup E_i$ denote the exceptional locus of the symplectic
contraction with $E_i$ being its irreducible components with the
reduced structure. In the present paper we shall discuss the case when
$E$ is contracted to an isolated singular point $y\in Y$. We will call
such a $\f$ {\sl a very small symplectic contraction}. The following
result can be found in \cite[Ch.~4]{Thesis} or
\cite[Sec.~4]{Kaledin00}, we sketch its proof after
\ref{thm-sigma-exact}.

\begin{lem}\label{exc-locus-dim} Let $\f: X\ra Y$ be a
very small symplectic contraction with exceptional locus $E=\bigcup_i
E_i$. Then the components $E_i$ are projective varieties of dimension
$\dim X/2$ and they are Lagrangian subvarieties of $X$.\end{lem}

By the Ionescu estimate, see \cite[IV.2.6]{Kollar}, in the above case
$E_i$ have the least possible dimension allowed for the components of
the exceptional locus of a symplectic contraction and therefore we
call it {\sl very small contraction}. If $\dim X=4$ then a very small
contraction is the same as a small contraction in terms of the Minimal
Model Program. By $\mu_i: \hat E_i \ra E_i\subset X$ we will denote
its normalization.

Although we are primarily interested in small symplectic contractions of
4-folds some of our preliminary arguments will be valid in arbitrary
dimension. We shall examine the structure of $E$ using two essential
tools: deformation of rational curves and vanishing theorems.


\subsection{Rational curves}

By Mori \cite{Mori79} and Kawamata \cite{Kawamata91}, the exceptional
locus of a crepant contraction is covered by rational curves. By
$\Hom(\P^1,X)$ we denote the scheme parametrizing morphisms from
$\P^1$ into $X$, its points are named after morphisms they are
associated to. For $f\in \Hom(\P^1,X)$ we can estimate the dimension
of all components of the parametrizing scheme from below by $\dim X +
\deg(f^*T_X)$, see \cite[II.1.2]{Kollar}. However, as it was observed
by Bogomolov, for symplectic manifolds one gets actually a better
estimate, see \cite{Ran} and \cite{Thesis}.

\begin{prop}\label{hom-dimension} Let $X$ be a 
symplectic manifold of dimension $n$. 
Suppose that either $X$ is projective or it admits 
a symplectic contraction onto an affine variety. 
Then the dimension of every
component of the scheme $\Hom(\P^1,X)$ is equal to $n+1$ at
least.\end{prop}

\noindent Let us sketch an argument, for details we refer to
\cite[Ch.~3]{Thesis} or to \cite[Sec.~2]{Wierzba}.  We consider a
1-parameter deformation of $X$ which does not extend the class of
$f(\P^1)$.  That is, in the compact case by a result of Bogomolov,
Tian and Todorov, while in the non-compact case by
\cite[Thm.~1.3]{Kaledin00}, see also \cite[Thm.~3.6]{KV}, there exists
a 1-parameter family of deformations of $X$ with the total space
$\ccX=\bigcup_t \ccX_t$ such that given $f:\P^1\ra X=\ccX_0$ no
deformation of it extends to $\ccX_t$ for $t\ne 0$. Hence we have the
equality $\Hom(\P^1,X)=\Hom(\P^1,\ccX)$ and thus the general estimate
applied to $\Hom(\P^1,\ccX)$ gives the above estime for
$\Hom(\P^1,X)$. \medskip

Let $F: \Hom(\P^1,X)\times\P^1\ra X$ denote the evaluation morphism
$F(f,p)=f(p)$. Consider an $Aut(\P^1)$-invariant subset
$V\subset\Hom(\P^1,X)$. By $F_V$ we denote the restriction of the
evaluation to $V\times\P^1$ and by $Locus(V)$ the closure of the image
of $F_V$. If $x\in F_V(V)$ then by $V_x$ we denote $V\cap
\Hom(\P^1,X;0\ra x)$ and by $Locus(V_x)$ the closure of $F(V_x)$. We
shall use the notions of an {\sl unsplit} and {\sl generically unsplit
subset} of $\Hom(\P^1,X)$, as defined in \cite[IV.2]{Kollar} (in case
of a generically unsplit irreducible component of $\Hom(\P^1,X)$ we
will consider its open and non-empty subset parametrizing birationally
rational curves). Moreover the term unsplit family will also be used
for a subset of $\Hom(\P^1,X;0\mapsto x)$ with the obvious change that
the subset has to be $Aut(\P^1;0)$-invariant.

\begin{defn}\label{minimal-dominating}
Let $X$ be a projective variety with an ample divisor $H$
and let $W\subset\Hom(\P^1,X)$ be an irreducible subset.
\item{(i)} By $\deg_H W$ we denote the degree of the $f^*H$
over $\P^1$ for a $f\in W$.
\item{(ii)} We say that $W$ dominates a subset $E\subset X$ if
$\deg_HW>0$ and $Locus(W)=E$.
\item{(iii)} We call $W$ a minimal dominating for $E$ if it satisfies (ii)
and $\deg_H W$ is minimal among all components of $\Hom(\P^1,X)$
which dominate $E$.
\end{defn}

We note that a minimal dominating $V$ is generically unsplit in the
sense of \cite[IV.2]{Kollar}. However, these two notions are not
equivalent. In fact, if $W$ is minimal in the above sense then for a
general $x\in Locus(V)$ the set $W_x=W\cap\Hom(\P^1,X;0\mapsto x)$ is
unsplit; such $W$ may be called {\sl locally unsplit}. If $X$ is
obtained by blowing up $\P^2$ at a point then the lift-up of a
(parametrization of a) general line produces a generically unsplit family
which is not locally unsplit in the above sense.

\begin{lem}\label{dominating-components} Let $\f: X\ra Y$ be a very
small symplectic contraction a with the exceptional locus
$E=\bigcup_iE_i$. Suppose that an irreducible component
$W\subset\Hom(\P^1,X)$ parametrizes curves contracted by $\f$ and
assume moreover that either $W$ is generically unsplit or $\dim X=4$
and $W$ contains a morphism which is birational onto its image. Then
there exists a (unique) irreducible component $E_W$ of $E$ which is
dominated by $W$
\end{lem}

\begin{pf} The proof follows from the above dimension estimate: we are
supposed to get $\dim Locus(V)\geq\dim X/2$. If $\dim X=4$ then we just
observe that $\dim V\geq \dim Aut(\P^1)$ so $\dim Locus(V)>1$. In
arbitrary dimension we we use \cite[IV.2.6.2]{Kollar}.\end{pf}

We recall that morphisms can be pushed-forward naturally. That is
given a morphism $\mu: M\ra N$ we have $\mu_*: \Hom(\P^1,M)\ra
\Hom(\P^1,N)$ defined as $\mu_*(f)=\mu\circ f$. If $\mu$ is the
normalization then for dominating subsets of $\Hom(\P^1,N)$ the
push-forward operation can be (generically) inverted.

\begin{lem}\label{lift-up} Let $\mu: \hat M \ra M$ be the normalization
of a projective irreducible variety and let $W\subset \Hom(\P^1,M)$ be
an irreducible subset which dominates $M$. Then there exists a
naturally defined subset $\hat W$ of $\Hom(\P^1,\hat M)$ which
dominates $\hat M$ and $\mu_*: \hat W\ra W$ is proper and
surjective. Moreover $\mu_*$ and the operation $W\mapsto \hat W$
define bijection between dominating components and $W$ is minimal
dominating component for $M$ if and only if $\hat W$ is minimal
dominating component for $\hat M$. \end{lem}

\begin{pf} Let $\mu_W: \tilde W\ra W$ be the normalization of $W$. Since
$W$ dominates $M$ the evaluation $F_W: W\times\P^1\ra M$ can be lifted
up to $\tilde F_W: \tilde W\times \P^1\ra \hat M$, so that
$\mu\circ\tilde F_W= F_W\circ(\mu_W\times id_{\P^1})$. Thus by the
universality of $\Hom$ we have a morphism $\psi: \tilde W\ra
\Hom(\P^1,\hat M)$ such that $\psi\times id_{\P^1}$ factors $\tilde
F_W$. We define $\hat W$ as the image of $\psi$ and because we have
$\mu_W=\mu_*\circ\psi$ then the properties of $\hat W$ follow.
\end{pf}

The above procedure will be called a lift-up to the normalization and will
be applied to normalizations $\mu_i: \hat E_i\ra E_i$ of components of
$E$. Similarly, by considering curves passing through a fixed point we
obtain families dominating normalizations and passing through a
fixed smooth point. Thus we get the following.

\begin{cor}\label{normalization-rational} Let $\f: X\ra Y$ be a very
small symplectic contraction with the exceptional locus
$E=\bigcup_iE_i$. Then the normalization of each irreducible component
$\hat{E_i}\ra E_i$ is rationally connected and it is a rational
surface if $\dim X=4$. \end{cor}

\begin{pf} Use \cite[IV.2.6.1]{Kollar} to get $\dim Locus(V_x)=\dim X/2$ for
$V$ generically unsplit and a general $x\in E_i= Locus(V)$.  Thus we
can lift up $V_x$ to a family of curves dominating $\hat E_i$.\end{pf}

In the course of the present paper we shall need the following
defintion, see \cite[Sect.3]{FOV}.

\begin{defn}\label{def-d-connected}
Let $Z$ be a variety or an analytic space with a closed point $z$. The
space $Z$ is called $d$-connected at $z$ if any component of $Z$
containing $z$ is of dimension bigger than $d$ and for any subspace
$T\subset Z$ containing $z$ of dimension smaller than $d$ there exists
a small neighborhood $U$ of $z$ in $Z$ such that $U\setminus T$ is
connected.\end{defn}

Now the following estimate on $d$-connectedness of $\Hom$ is obtained
exactly as the dimension estimate.

\begin{prop}\label{hom-d-connected} Let $X$ be a projective variety
which contains a rational curve $f:\P^1\ra C\subset X$. Then the scheme
(or the analytic space) $\Hom(\P^1,X)$ is $\dim X + \deg(f^*TX) - 1$
connected at $f$.\end{prop}

\begin{pf} By \cite[I.2.8 and proof of I.2.16]{Kollar} the $\Hom$ scheme is
defined locally by $\dim H^1(\P^1,f^*TX)$ equations in a smooth space of
dimension $\dim H^0(\P^1,f^*TX)$ thus the proposition follows by
\cite[3.1.13]{FOV} and Riemann-Roch. \end{pf}

In the case of a symplectic manifold the argument of extending $X$ by
taking its 1-parameter deformation (which we presented above arguing for
\ref{hom-dimension}) applies and the actual estimate is better. 

\begin{cor}\label{hom-4-connected} Let  $X$ be  a symplectic $2n$-fold. 
Suppose that either $X$ is projective or it admits a symplectic 
contraction onto an affine variety.
Then $\Hom(\P^1,X)$ is $n$-connected at any point. \end{cor}


\subsection{Vanishings}

We need the following vanishing due to Grauert-Riemenschneider, Kodaira,
Kawamata and Viehweg, see \cite[Sect.~1-2]{KMM}.

\begin{thm}\label{GRKKV-vanishing} Let $X$ be a smooth variety and
$\f:X\ra Y$ a birational proper morphism. If $L$ is a line bundle such
that $K_X+L$ is $\f$-big and nef then $R^i_*\f L=0$ for $i>0$. In
particular if $\f$ is a Fano-Mori or crepant contraction then
$R^i\f_*\ccO_X=0$ for $i>0$.\end{thm}

The vanishing is needed, among other things, for the following.

\begin{lem}\label{picard=H2} Let $\f: X\ra Y$ be a Fano-Mori or crepant
birational contraction of a smooth variety with the exceptional locus
$E=\bigcup_i E_i$ which is contracted to a point $y\in Y$. Then after
possible shrinking both $X$ and $Y$ to an analytic neighborhood of $E$
and $y$, respectively, we have $Pic(X)\iso H^2(E,\Z)$.
\end{lem}

\begin{pf} Firstly, because of the above vanishing, we may shrink $X$ so
that $H^i(X,\ccO_X)=0$, for $i>0$, and therefore the Chern class map
$Pic(X)\ra H^2(X,\Z)$ becomes an isomorphism. Secondly, again possibly
shrinking $X$, we get $E$ a deformation retract of $X$, by
\cite{Lojasiewicz}, and therefore $H^i(E,\Z)\iso H^i(X,\Z)$. Combining
these two we are done.
\end{pf}

We shall need the following.

\begin{prop}\label{top-vanishing} Let $\f: X\ra Y$ be a very small
contraction of a symplectic manifold of dimension $2n$. Let  $E=\bigcup_iE_i$ be
the exceptional locus. By  $\mu_i:\hat{E_i}\ra E_i$ we denote the
normalization of irreducible components. Then
$$H^n(E_i,\ccO_{E_i})=H^n(\hat{E_i},\ccO_{\hat{E_i}})=0$$ \end{prop}

\begin{pf} The vanishing for $E_i$ is in \cite[1.7]{AW-SC}. 
Then the vanishing for
the normalization follows by cohomology of the sequence
$0\ra\ccO_{E_i}\ra(\mu_i)_*\ccO_{\hat{E_i}}\ra\ccQ\ra 0$, where the
support of $\ccQ$ is of dimension smaller than $n$. \end{pf}

\begin{cor}\label{normalization-sing} Let $\f: X\ra Y$ be a small
contraction of a symplectic 4-fold with the exceptional locus
$E=\bigcup_iE_i$. Then the normalization of each irreducible component
is a rational surface with quotient singularities.
\end{cor}

\begin{pf} Let $\pi_i: \tilde{E_i}\ra\hat{E_i}$ be a desingularisation.
Then, by rationality \ref{normalization-rational} we have
$H^j(\tilde{E_i},\ccO_{\tilde{E_i}})=0$ for $j=1,\ 2$ and by
\ref{top-vanishing} we get additionally
$H^2(\hat{E_i},\ccO_{\hat{E_i}})=0$ so $R^1_*\pi_i\ccO_{\tilde{E_i}}=0$
by Leray spectral sequence for $\pi_i$. Thus $E_i$ has rational
singularities which are quotient singularities as well.
\end{pf}


\subsection{Normal surfaces with many rational curves}

The main result of this section is a 2-dimensional version of a result
by Cho, Miyaoka and Shepherd-Barron announced in earlier versions of
\cite{CMSB}. We include the proof to make the present paper
self-contained. We begin by stating a characterization of $\P^2$
obtained in the course of the proof of the main theorem of
\cite{Kawamata89}, see step 2.2 of the proof.

\begin{prop}\label{kawamata-P2} Let $S$ be a normal projective surface
with an ample divisor $H$. Let $f: \P^1\ra S$ be a morphism whose
degree $\deg(f^*H)$ with respect to $H$ is minimal but positive. If
for some $p\in\P^1$, with $f(p)\in S\setminus Sing(S)$, we have
$\dim_{f}\Hom(\P^1,S; p\mapsto f(p))\geq 3$ then
$S\iso\P^2$.\end{prop}

We need to extend this characterization to quotients of $\P^2$.

\begin{thm}\label{cmsb} Let $S$ be a normal projective surface with
quotient singularities. Suppose that for a smooth point $s\in S$ there
exists a component $W_s$ of $\Hom(\P^1,S; 0\mapsto s)$ which is unsplit and
$\dim W_s\geq 3$.
Then the fundamental group $G=\pi_1(S\setminus Sing(S))$ acts
algebraically on $\P^2$ and $S$ is the quotient $\P^2/G$ with the quotient
morphism $\nu: \P^2 \rightarrow S$ which is smooth covering
 outside the inverse
image of singularities of $S$. \end{thm}

In the course of the proof we will construct the morphism $\nu$. We note
that once the covering $\nu: \P^2 \rightarrow S$, which is smooth outside
$\nu^{-1}(Sing(S))$, is constructed then the rest of the theorem follows.
Indeed, its restriction $\P^2\setminus\nu^{-1}(Sing(S))\ra S\setminus
Sing(S)$ is then the universal cover of the smooth locus of $S$. Thus, if
$G$ is the Galois group of this covering then it acts on $\P^2$ and 
$S=\P^2/G$.

We follow Koll\'ar's book \cite[pp. 108--112]{Kollar} and construct a
normal family of rational curves $V$ (we drop Koll\'ar's subscript
$^n$) which is obtained from the normalization of $W_s$ as a quotient
by the action of $Aut(\P^1,0)$ with the universal $\P^1$ bundle
$\pi_U: U\rightarrow V$ which admits the evaluation morphism $F_V:
U\rightarrow S$.  Points of $V$ will be denoted by classes of
morphisms, that is by $[f]$, where $f\in W_s$.  The bundle
$U\rightarrow V$ has a section $V_0\subset U$ which is contracted by
$F_V$ to $s$.  Since the family $W_s$ is unsplit the quotient $V$ is
proper and the morphism $F$ is finite-to-one outside $V_0$ hence it is
surjective.

After having made this preliminary construction, to which we will refer
in the course of the proof, we state a somewhat more general observation.

\begin{lem} Let $B$ be a germ of a smooth curve (a small disc in the analytic
set-up) with the closed (central) point $b$.  Consider a morphism $F:
B \times \P^1 \rightarrow X$ into a smooth variety such that
$dim(imF)=2$.  By $f$ let us denote $F_{|b\times\P^1}$. Let $p_1,
\dots, p_k \in\P^1$ be such that $F(B \times {p_i}) = x_i$. Then the
quotient $(f^*T_X)/T_{\P^1}$ has a (generically) non-zero section
vanishing at points $p_i$.
\end{lem}

\begin{pf}
Let $t$ be a local coordinate in $B$ with $b=\{t=0\}$ and $v$ a
nonvanishing vector field.  Via the tangent morphism $T_{B \times \P^1} =
T_B \times T_{\P^1}\rightarrow F^*(T_X)$ the field $v$ gives a
non-zero section $s$ in $F^*(T_X)/T_{\P^1}$.  The section vanishes
along $B \times {p_i}$. Let $m$ denote the multiplicity with which it
vanishes along $b \times \P^1$ ($m$ may be zero as well).  Then
$s/t^m$ does not vanish identically along $b \times \P^1$ and it
is zero along $B \times {p_i}$ so its restriction to $b \times \P^1$
is what we are looking for.
\end{pf}

\begin{cor}\label{smoothness-over-smooth}
In the notation introduced above let $f\in W_s$ be such that
$f(\P^1)$ is contained in the smooth locus of $S$. Then $f^*(T_S) =
\ccO(2) \oplus \ccO(1)$ and consequently $f$ is immersion
(that is $f(\P^1)$ has no cusp but possibly nodes), $W_s$ is smooth at
$f$, and the evaluation $F_{W_s}$ is of maximal rank (smooth) along $f
\times (\P^1 -\{0\})$. In particular $V$  is smooth at $[f]$, and 
$F_V: U\rightarrow S$ is \`etale over $\pi_U^{-1}([f])
\setminus V_0$.
\end{cor}

\begin{pf}
Since $W_s$ is unsplit $deg(f^*(T_S))\leq 3$. Next, by
the previous lemma, $f^*(T_S) = \ccO(2) \oplus \ccO(1)$. The rest
follows by \cite[II.3.4]{Kollar} and
\cite[II.2.16]{Kollar}.
\end{pf}

Now we want to have a version of the above corollary also to $f$ such that
$f(\P^1)$ meets singularities of $S$. For this purpose we need a general,
somewhat technical, observation which holds in the analytic category. 

\begin{lem}\label{technical-covering}
Let $X$ be a normal complex variety with isolated quotient 
singularities. Let $C=\P^1$ and
let $\Delta$ denote a small disc around $0 \in \C$. Let $F: \Delta \times C
\rightarrow X$ be a holomorphic morphism such that $F(\Delta\times C
\setminus \bigcup_i \{(0,p_i)\}$ is contained in the smooth locus
of $X$. By $f$ we denote $F_{|(0\times C)}$ and by $x_i$ we denote
$F(0,p_i)$. Suppose that $f$ is birational onto its image.  Then,
after possibly shrinking $\Delta$ to a smaller disc $\Delta'\ni 0$ we can
choose an open (analytic) subset $\ccU\subset X$ such that $F(\Delta'\times
C)\subset \ccU$ and there exists a complex manifold $\ccV$ and a
holomorphic morphism $\pi: \ccV\rightarrow \ccU$ which satisfies the
following conditions:
\item{(i)} $\pi$ is finite and smooth outside $\bigcup_i \pi^{-1}(0,x_i)$,
\item{(ii)} there exists a morphism $F': \Delta'\times C \rightarrow \ccV$
such that $\pi\circ F' = F_{|\Delta' \times C}$.
\end{lem}

\begin{pf}
First we claim that we can find an \'etale cover $\ccV' \rightarrow
\ccU$ of a neighborhood of $f(C)$ so that after lifting $F$ to $F':
\Delta'\times C\ra \ccV'$ the resulting $f': C \rightarrow \ccU'$ is
bijective onto its image.  In other words we want to separate branches
of nodes of $f(C)$. (We can make our argument assuming that $\ccU$ is
smooth since at this stage we can embed all data into a smooth
variety.) So we take $\ccU$ to be a small tubular neighborhood of
$f(C)$ such that if $x=f(p_i)$, with $i=1,..k$, is a multiple point
then there exists a neighborhood $\ccU_x$ (say a ball) such that $\ccU
\cap \ccU_x = \ccU_x^1\cup...\cup \ccU_x^k$ and each $\ccU_x^i$ is a
tubular neighborhood of the respective branch of $C$. That is
$C^i=f^{-1}(\ccU_x^i)$ are disjoint neighborhoods of $p_i$, for
$i=1,\dots, k$, and each $\ccU^i_x$ is a tubular neighborhood of
$f(C^i)$. We can moreover assume that the intersection of the boundary
of the ball $\ccU_x$ with $\ccU$ splits into disjoint connected
components each one of them being the intersection of the boundary of
$\ccU_x$ with $\ccU_x^i$.  Now around $x$ the new variety $\ccU'$ can
be obtained by taking the disjoint union of $\ccU_x^i$, the glueing of
the components to the rest of $\ccU$ is obvious. It is also clear how
to lift up $f$ and $F$ to $f': C \rightarrow \ccU'$ and $F':
\Delta\times C\rightarrow \ccU'$, respectively. We repeat this
construction for any multiple image point $x\in f(C)$, and to
simplify the notation, we call the result $\ccU$ again etc.

Thus, from now on we may assume that $f$ is bijective onto its image.
Let $p\in C$ be a point mapped to a singular $x\in \ccU$. Let $\ccU_x$
be a neighborhood of $x$ such that there exists a morphism $\pi_x:
\ccV_x \rightarrow \ccU_x$ with $\ccV_x$ smooth and $\pi_x$ etale
outside $y=\pi_x^{-1}(x)$. Now, by possibly shrinking $\ccU_x$, we can
cover $f(C)$ by $\ccU_x$ and an open set $\ccV'\subset X$ not
containing $x$ which is a tubular neighborhood of $f(C)$ at its
intersection with $\ccU_x$. That is we want $\ccU_x\cap \ccV'= f(A)
\times D$ where $A$ is an annulus around $p$ and $D$ is
contractible. Let $D_1$ be a small disc around $p$ such that $F$ maps
$\Delta\times D_1$ into $\ccU_x$. Thus we have $\Delta\times
D_1\setminus (0,p)$ mapped into $\ccU_x\setminus\{ x\}$ and since it
is simply connected the map lifts up to the covering space
$\ccV_x\setminus\{y\}$. This implies that the loop generating
the fundamental group of $\ccU_x\cap \ccV'= f(A) \times D$ lifts up to
$\ccV_x\setminus \{y\}$ hence the inclusion $\ccU_x\cap
\ccV'\hookrightarrow \ccU_x\setminus\{ x\}$ lifts up to $\ccU_x\cap
\ccV'\rightarrow \ccV_x\setminus\{ y\}$. Now we define $\ccV$ as
$\ccV'$ attached to $\ccV_x$ via this lift-up. The definition of $F':
\Delta'\times C \rightarrow \ccV$ comes with the construction.  Now we
can apply this construction consecutively to all singular points $x_i$
which lie on $f(C)$ and we conclude the lemma.
\end{pf}

Now we continue with the proof of \ref{cmsb}. We now extend 
\ref{smoothness-over-smooth} to the following.

\begin{lem} The morphism $F: U\rightarrow S$ is \'etale outside 
$F^{-1}(s)\cup F^{-1}(SingS)$.\end{lem}

\begin{pf} The smoothness of the morphisms $F$ along fibers of $\pi_U$
which are mapped into the smooth locus of $S$ was already observed in
\ref{smoothness-over-smooth}. The argument leading to it  is actually local
and holds in analytic category too. Thus it can be applied to the local
lift-up family of rational curves in $\ccV$ which we have constructed in
the preceding lemma.  That is, in the notation of the previous lemma, if
$f'=F'_{|(0\times \P^1)}: \P^1 \rightarrow \ccV$ is the liftup of $f$
then $deg(f')^*(T\ccV)=3$ and, by \ref{smoothness-over-smooth}, 
$(f')^*(T\ccV)=\ccO(2)\oplus\ccO(1)$. Consequently the rest of the
argument follows.  \end{pf}

\noindent {\it Conclusion of the proof of \ref{cmsb}.}  Let $U\rightarrow
S'\rightarrow S$ be the Stein factorization of $F$, the morphism
$U\rightarrow S'$ contracts $V_0$ to a point which is the vertex of the
cone.  By the previous lemma $S'\rightarrow S$ is etale in codimension 1.
Since $s$ is a smooth point then by purity of the branch locus it
follows that $S'\rightarrow S$ is etale over $s$ hence $V_0$ is
contracted to a smooth point on $S'$. But $S'$ is a cone with the vertex
of the image of $V^0$ hence it must be $\P^2$ (and $U$ is the 1st
Hirzebruch surface).  Thus we have $S' = \P^2 \rightarrow S$  \'etale in
codimension 1.

\bigskip

We restate \ref{cmsb} to the following.

\begin{thm}\label{cmsb1} Let $S$ be a normal projective surface with
quotient singularities and $H$ an ample divisor over $S$. Assume that
$W\subset \Hom(\P^1,S)$ is a minimal dominating component for $S$. If
$\dim W\geq 5$ then $S\iso\P^2/\pi_1(S\setminus Sing(S))$  is as
described in the conclusion of \ref{cmsb} and moreover the quotient map
$\nu: \P^2\ra S$ induces a surjective  morphism of components of
$\Hom$-schemes $\nu_*: W_1\ra W$, where  $W_1\subset\Hom(\P^1,\P^2)$
parametrizes lines.  In particular such $W$ is unique. If moreover {\bf any}
$f\in W$ is birational  onto the image then $S\iso\P^2$. \end{thm}

\begin{pf} By what we have noticed above $W_s=W\cap \Hom(\P^1,S;
0\mapsto s)$ is unsplit for a general choice of $s\in S$ and of
dimension $\geq 3$. Thus \ref{cmsb} applies and actually its proof
produces lifting of curves parametrized by $W_s$ to lines on $\P^2$.
Thus it remains to prove the last statement of the theorem. To this end
we note that the action of $G=\pi_1(S\setminus Sing(S))$ on $\P^2$
induces a dual action on $(\P^2)^*$ which parametrizes lines on $\P^2$.
(The quotient $(\P^2)^*/G$   can be identified in $\Chow(S)$ as a
component  parametrizing images of lines and the morphism $\nu_*:
(\P^2)^*\ra (\P^2)^*/G$ is the push-forward of cycles.) In particular,
if $G_l\subset G$ is  the stabilisor of a line $l\in (\P^2)^*$ then the
morphism $\nu_{|l}$  is of a degree equal to $|G_l|$. Thus,
if $\nu_l$ is birational then  $(\P^2)^*\ra (\P^2)^*/G$ is unramified covering
hence it is isomorphism and we are done.

\end{pf}

Having completed the proof of \ref{cmsb} we derive the following
conclusion for the proof of our main theorem.

\begin{cor}\label{normalization-surface} Let $\f: X\ra Y$ be a small
contraction of a symplectic 4-fold with the exceptional locus
$E=\bigcup_iE_i$. Then the normalization of each irreducible component
$\mu_i: \hat{E_i}\ra E_i\subset X$ admits a finite morphism
$\nu_i:\P^2\ra\hat{E_i}$ which has the properties described in \ref{cmsb}
and \ref{cmsb1}.
\end{cor}

\begin{pf} In view of \ref{hom-dimension}, \ref{dominating-components}, 
\ref{lift-up} and
\ref{cmsb} the description of normalization follows.  \end{pf}


\subsection{A vector field}

We begin this section by discussing a general fact related to
properties of Fano-Mori and crepant contractions. This property was
also observed by Campana and Flenner \cite{Campana-Flenner}.

\begin{prop}\label{thm-sigma-exact} Let $\f: X\ra Y\ni y$ be a local
analytic Fano-Mori or crepant contraction and let $\sigma\in \Omega^2_X$
be a closed holomorphic form. After possibly shrinking $X$ to a smaller
neighborhood of $\f^{-1}(y)$ there exists a $1$-form  $\alpha\in
H^0(X,\Omega_X)$, such that $d\alpha=\sigma$. \end{prop}

\begin{pf}

We consider the analytic de Rham complex over $X$: $$0\ra \C_X \ra
\ccO_X \stackrel{d^0}{\lra} \Omega^1_X \stackrel{d^1}{\lra}   \Omega^2_X
\stackrel{d^2}{\lra}\Omega^3_X \stackrel{d^3}{\lra}\cdots$$ Let us set
$\ccK=\im(d_0)=\ker(d_1)$.  In order to prove the proposition
we have to show that $H^1(X,\ccK)\ra H^1(X,\Omega^1_X)$ is injective.
After possibly shrinking $X$ so that its image is Stein we have vanishing
$H^i(X,\ccO_X)=0$ for $i>0$. This implies that the boundary map
$\delta: H^1(X,\ccK)\ra H^2(X,\C)$ is an isomorphism.

Now we look at the exponential sequence $0\ra \Z_X\ra \ccO_X\ra
\ccO^*_X\ra 0$ to find out, again by  $H^i(X,\ccO_X)=0$ for $i>0$,
that the boundary $\delta': H^1(X,\ccO^*_X)=H^2(X,\Z_X)$
is an isomorphism too. In fact $\delta$ and $\delta'$ are related
by a commutative square appearing in the following diagram 
\begin{eqnarray*}
\begin{CD}
\H^1(X,\ccO^*_X)@>>> \H^1(X,\ccK) @>>> \H^1(X,\Omega^1_X)\\
@V{\delta'}VV @V{\delta}VV \\
\H^2(X,\Z)@>{\otimes\C}>>\H^2(X,\C).
\end{CD}
\end{eqnarray*}
where the composition of the upper row arrows is the Chern
class $c: \H^1(X,\ccO^*_X)=\Pic X\ra H^1(X,\Omega^1_X)$, with
$c(g)=dg/g=d(log(g))$, see
\cite[II Ex.~7.4]{Hartshorne}. Thus we will be done if $c$ is injective.

Suppose $L\in\Pic X$ is a line bundle such that $c(L)=0$. The
definition of $c$ is functorial hence descends to any curve contracted
by $\f$ so $L$ is numerically trivial. This, by Kawamata-Shokurov base
point freeness theorem, see \cite[3-1-1, 3-1-2]{KMM} or
\cite[3.24]{KM}, implies that $L$ is actually trivial. \end{pf}

At this point let us sketch a proof of the dimension 
estimate on the components
of the exceptional locus of a symplectic contraction \ref{exc-locus-dim}, 
see \cite[Ch. 4]{Thesis}. 

\noindent{\it Proof of \ref{exc-locus-dim}.}  The proof comes by
combining the two tools which we have explained above.  From the
dimension estimate on $\Hom$, \ref{hom-dimension}, we get Ionescu's
inequality, see \cite[IV.2.6]{Kollar}: $dimE_i\geq n$. On the other
hand by \ref{thm-sigma-exact}, after shrinking $X$ to a small
neighborhood of $E$ we see that the cohomology class of $\sigma$ in de
Rham cohomology of $X$ is trivial.  Let $\tilde{E_i}\ra E_i$ be a
projective resolution of singularities of a component $E_i\subset E$.
The pullback of the form $\sigma$ is zero in the cohomology of $\tilde{E_i}$
hence, by Hodge theory on $\tilde{E_i}$, the pull-back is zero 2-form
itself.  Thus for any smooth point $x\in E_i$ the linear space
$T_xE_i\subset T_xE$ is isotropic for $\sigma$, hence
\ref{exc-locus-dim}.

\medskip

The \ref{thm-sigma-exact} will work together with the following results
which we formulate in a somewhat more general context.

Let $X$ be a symplectic $2n$-fold with the symplectic form $\sigma$.
Suppose that there is a $1$-form $\alpha$ such that $d\alpha=\sigma$.
Let $\xi$ be the vector field on $X$ which is via $\sigma$ dual to
$\alpha$. Let $M\subset X$ be a  (reduced, but possibly reducible)
Lagrangian subvariety such that for some resolution $M'\ra M$ we have
$H^0(M',\Omega_{M'})=0$. Then the following two results hold.

\begin{prop}\label{thm-xi-pres-e}
The vector field $\xi$ preserves $M$, that is, it induces a derivation 
$\zeta\in \H^0(M,\Theta_M)$.
\end{prop}

\begin{pf}
Let $\ccI\subset \ccO_X$ be the ideal sheaf defining $M$. Since $M$ is
Lagrangian, there is a commutative diagram
\begin{eqnarray*}
\begin{array}{ccccccc}
     & \ccI/\ccI^2 & \ra & \Omega_X|M    & \stackrel{a}{\ra} & \Omega_M & \ra 0 \\
     & \da         &     &  ||\hat\sigma &     &  \da \\
0\ra & \Theta_M    & \ra & T_X|M      & \stackrel{b}{\ra}& (\ccI/\ccI^2)^{\vee} \\
\end{array}
\end{eqnarray*}
where $\xi=\hat\sigma(\alpha)$.
Suppose $a(\alpha)$ is nonzero at some generic point
of $M$. Then $a(\alpha)$ induces a nonzero section of $\Omega_{M'}$. 
But by assumption $H^0(M',\Omega_{M'})=0$. 
Therefore $a(\alpha)$ is zero at all generic points of $M$. Since
$(\ccI/\ccI^2)^{\vee}$ is torsion free it follows that $b(\xi)$ is zero
and therefore $\xi$ lifts to $\zeta\in H^0(M,\Theta_M)$.
\end{pf}

The following result is due to Dmitri Kaledin.

\begin{prop}\label{thm-kal}
In the above situation if the derivation $\zeta$ is zero on $M$ then $M$ is smooth.
\end{prop}

\begin{pf}
$M$ is Lagrangian, hence of dimension $n$. If $\zeta$ is zero on $M$,
then the image of $\alpha$ in $\Omega_X|M$ is zero and therefore
$\alpha$ lifts to a global section of $\Omega_X\tensor \ccI$, where
$\ccI\subset \ccO_X$ is the ideal sheaf defining $M$.  Let $p\in M$ be
a point. Then $\alpha_p\in \Omega_{X,p}$ is of the form $\alpha=\sum
a_i db_i$ with $a_i\in \ccI_p$ and $b_i\in \ccO_{X,p}$. Since $\sigma$
is symplectic,
$$\sigma^{\wedge n} =(d\alpha)^{\wedge n} 
  =\sum da_{i_1}\wedge db_{i_1}\wedge\cdots \wedge da_{i_n}\wedge db_{i_n}$$
is a non vanishing section of $\omega_X$. Let $\m_p\subset \ccO_{X,p}$
be the maximal ideal. If $M$ is non-smooth at $p$, then
the image of $\ccI_p\ra \m_p/\m_p^2=\Omega_X\tensor k(p)$ is at most
$(n-1)$-dimensional. Therefore 
$da_{i_1}\wedge\cdots\wedge da_{i_n} \mod \m_p=0$ 
for all $i_1,\dots,i_n$.
Therefore $\sigma^{\wedge n}$ vanishes at $p$. Contradiction.
\end{pf}

As an application of the above we prove the following key ingredient
to the proof of the main theorem.

\begin{prop}\label{normalization-final} Let $\f: X\ra Y$ be a
small contraction of a symplectic 4-fold with the exceptional locus
$E=\bigcup_iE_i$. Let $Sing(E)$ denote the singular locus of $E$. Then
the normalization of each irreducible component $\mu_i:\hat{E_i}\ra E_i$
admits a uniquely defined finite morphism $\nu_i:\P^2\ra\hat{E_i}$ which is
unramified in codimension 1. Moreover  the inverse image of $Sing(E)$ 
under the composition $\psi_i:\P^2\ra\hat{E_i}\ra E_i\subset E$
consists of points and rational curves. \end{prop}

\begin{pf} In view of \ref{normalization-surface} we are only to prove
the statement about the inverse image of the singular set. By
\ref{thm-sigma-exact} we can construct on $X$ a 1-form $\alpha$ such
that $d\alpha=\sigma$.

Suppose first that $\alpha$ does not vanish identically on the
components $E_i$. Then by \ref{thm-xi-pres-e} it produces non-trivial
differentiation $\zeta$ on $E_i$. The differentiation $\zeta$ lifts up
to $\hat{E_i}$ and because the morphism $\P^2\ra\hat{E_i}$ is \'etale
outside the inverse image of singularities of $\hat{E_i}$ it litfs up to
$\P^2$ too. Another explanation is as follows: by the uniqueness of the
construction of the morphism $\P^2\ra\hat{E_i}\ra E_i\subset E\subset X$
the action of a 1-parameter group on $X$ associated to the vector field
$\xi$ lifts up to $\P^2$. Now by its very construction the
differentiation (or the 1-parameter group action) must preserve the
inverse image of the singular set of $E$. Thus we are only to note that
among curves on $\P^2$ only rational have this feature.

Now suppose that $\alpha$ vanishes on a component $E_i$ then by
\ref{thm-kal} $E_i=\P^2$ and it can not be identically zero on any
components of $E$ meeting $E_i$. Thus by the previous part the common
locus consist of isolated points or rational curves. \end{pf}


\section{Proof of Theorem~\ref{big-thm}}

From now in this section on we deal exclusively with small contractions
of symplectic 4-folds. That is, by $\f:X\ra Y$ we will denote a small
contraction of a symplectic 4-fold with a connected exceptional locus
$E=\bigcup_i E_i$ contracted to an isolated singular point $y\in Y$.

At this point, in view of \ref{normalization-final} we have some
information on the normalization of components $E_i$ as well as on the
locus of their common points and singularities. We want to prove that
actually there is only one component and it is $\P^2$. Our
intermediate task will be the irreducibility of $E$ and
\ref{normal-homeo} which says that the normalization of $E$ is a
homeomorphism to $\P^2$ in codimension 1. Because of \ref{picard=H2}
the line bundle $\ccO(1)$ from $\P^2$ extends over $X$ and Kawamata's
base point free technique can be applied. This idea was suggested to
us by Daniel Huybrechts and Manfred Lehn.


\subsection{Reduction to $\rho(X/Y)=1$}

Let $E=F_1\cup \cdots \cup F_r$ be the decomposition of $E$ 
into 1-connected components, that is each $F_i$ is 1-connected and
$F_i$ meets $F_j$ in a finite number of points , see \cite[3.1.5]{FOV}.

\begin{lem}\label{factorization} After possibly shrinking $X$, the rank
of $Pic(X)$ is $r$ and for any $i=1,\dots,r$ there is a contraction
$\f_i:X\ra Y_i$ which factorizes $\f:X\ra Y$ and such that the exceptional
set of $\f_i$ is precisely $F_i$. \end{lem}

\begin{pf} Firstly we note that because of \ref{normalization-surface}
all curves in $F_i$ are numerically equivalent, that is $\dim
N_1(F_i)=1$. Since $F_i$ meet in dimension 0 the natural surjection
$\coprod_i F_i\ra E$ implies an isomorphism $\bigoplus_i N_1(F_i)\iso
N_1(E)$ under which the cone $N\! E(E)$  of effective 1-cycles in $E$ is
generated by classes of effective cycles in $F_i$. Thus $N_1(E)$ is of
dimension $r$ and $N\! E(E)$ is simplicial with rays generated by
classes of curves in $F_i$. On the other hand by \ref{picard=H2}, after
possibly shrinking $X$, the inclusion $E\subset X$ implies the isomorphism
$N^1(X/Y) \iso N^1(E)$. Thus $\rho(X/Y)=r$ and the rest of the lemma
follows by the contraction theorem, see e.g. \cite[3.7]{KM}.\end{pf}

\begin{prop}\label{prop-red-rho-i} Reduction to $\rho(X/Y)=1$: Suppose
that Theorem~\ref{big-thm} holds with the additional assumption
$\rho(X/Y)=1$ (or with equivalent assumption that $E$ is 1-connected).
Then it is true without this assumption as well.
\end{prop}

\begin{pf} Let $f:X\ra Y$ be a contraction as in Theorem~\ref{big-thm}.
Suppose that $\rho(X/Y)=r>1$. Let $F_1,\dots,F_r$ be as above. Then by
\ref{factorization}we obtain for each $i=1,\dots,r$ a contraction
$f_i:X\ra Y_i$, which contracts just $F_i$. Note that $\rho(X/Y_i)=1$.
Therefore we can apply Theorem~\ref{big-thm} with the additional 
assumption $\rho=1$ to the contraction morphism $f_i:X\ra Y_i$ to
conclude that $F_i\cong \P^2$. 

Therefore all exceptional components of $E$ are $\P^2$'s and they meet
only in finitely many points. Let $E_1$ and $E_2$ be two irreducible
components with nonempty intersection. Let $X'$ be be obtained from
$X$ by the Mukai flop of $E_1$. Then the strict transform $E'_2$ of
$E_2$ is a nontrivial blow up of a $\P^2$. On the other hand there is
a contraction $X'\ra Y$ satisfying the conditions of
Theorem~\ref{big-thm} and $E'_2$ is an exceptional
component. Therefore, by \ref{normalization-surface} the surface
$E'_2$ cannot be a nontrivial blow up of $\P^2.$ Contradiction.
\end{pf}


\subsection{Connectivity arguments}\label{connectivity-section}

In view of Proposition~\ref{prop-red-rho-i} we assume from now  on that
$\rho(X/Y)=1$ which, by \ref{factorization}, is the same as saying that
$E$ is is connected in codimension  $1$. By $L$ let us denote the
positive generator of $Pic(X/Y)$. Let $E=E_1\cup \cdots\cup E_n$ be the
decomposition into  irreducible components. For each $i=1,\dots, n$, let
$\mu_i:\hat{E_i} \ra E_i$ be the normalization and let 
$\nu_i:\P^2\ra \P^2/G_i=\hat{E_i}$ be the quotient map given by 
\ref{normalization-surface} and let $\psi_i:\P^2\ra E_i$ be the 
composition of these two maps. By $(\nu_i)_*$ and $(\mu_i)_*$ we denote
the natural morphism of $\Hom$-schemes: $\nu_*(f)=\nu\circ f$,
$\mu_*(f)=\mu\circ f$.

The task of the present section, which is proving that $E$ is homeomorphic
to $\P^2$ in codimension1, is achieved in three steps. The main tool used
in each step is the connectivity feature \ref{hom-4-connected}.

\begin{prop}\label{normalization-P2} The normalization $\hat{E_i}$
of any component is isomorphic to $\P^2$.\end{prop}

\begin{pf} Let $W_i\subset\Hom(\P^1,X)$ be a minimal component
with respect to $L$ dominating $E_i$. Using the normalization we lift up
$W_i$ to $\hat W_i\subset \Hom(\P^1,\hat E_i)$ which is minimal and
dominating for $\hat E_i$ and thus by \ref{cmsb1} it is unique and the
image $(\nu_i)_*(W^1)$ of $W^1\subset \Hom(\P^1,\P^2)$ parametrizes lines.
Suppose that $\nu_i$ is not isomorphism. We will prove that this implies
that $\Hom(\P^1,X)$ is not 4-connected at some point hence we will
arrive to contradiction with \ref{hom-4-connected}.

According to \ref{cmsb1} there exists $f\in W^1$ such that the degree of the
morphism $\pi_i\circ f=\mu_i\circ\nu_i\circ f$ is bigger than 1, say it
is of degree $d$. Let $g\in\Hom(\P^1,X)$ denote the normalization of the
image of $\psi_i\circ f$. We have a degree $d$ covering $\pi_d:
\P^1\ra\P^1$ such that $g\circ \pi_d=\psi_i\circ f$. Let
$W_g\in\Hom(\P^1,X)$ be a component containing $g$, then by
\ref{dominating-components} there exists a component $E_j\ne E_i$ which
is dominated by $W_g$. Now taking composition $W_g\ni h \mapsto h\circ
\pi_d\in \Hom(\P^1,X)$ we obtain a subset of $\Hom(\P^1,X)$ which
dominates $E_j$ and which contains also $\psi_i\circ f\in W_i$. Thus,
apart of $W_i$ we have at least one irreducible component of
$\Hom(\P^1,X)$ which contains $\psi_i\circ f$.

Now in order to prove that $\Hom(\P^1,X)$ is not 4-connected at
$\psi_i\circ f$ it is enough to show that if $h$ is a small deformation
of $\psi_i\circ f$ which is in $W_i\setminus \psi_i\circ f\circ
Aut(\P^1)$ then $h$ is not contained in any other irreducible component
of $\Hom(\P^1,X)$. If however it was the case then we would get another
component of $\Hom(\P^1,X)$ which dominates $E_i$ and which is of the
same degree with respect to $H$ as the component $W_i$. This contradicts
the uniqueness statement of \ref{cmsb1}. 
\end{pf}

\begin{prop}\label{$E$-irreducible} $E$ is irreducible.\end{prop}

\begin{pf} Suppose that the proposition is not true. Since $E$ is
1-connected we can choose $C\subset E$, a curve which is common to
components $E_1,\dots,E_r$ where $r> 1$. Let $f_1: \P^1\ra C\subset X$
be the normalization. It is uniquely defined up to the action of
$Aut(\P^1)$.  Let $W_1$ be an irreducible component of $\Hom(\P^1,X)$
which contains $f_1$. Then, because of \ref{lift-up} the evaluation
map $F_{W_1}$ dominates one of the components of $E$, say $E_1$. Let
$E_2$ be another component of $E$ which contains $C$ as well. Let
$\mu_2: \P^2\ra E_2\subset X$ be the normalization,
\ref{normalization-P2}. Take a curve $C_2\subset\P^2$ which is mapped
to $C$. By \ref{normalization-final} the curve $C_2$ is rational so
let $f_2:\P^1\ra C_2$ be its normalization; by $\hat W_2\subset
\Hom(\P^1,\P^2)$ denote the (unique) component containing $f_2$. Let
$W_2=(\mu_2)_*(\hat W_2)$ be its image which dominates $E_2$. As in
the proof of \ref{normalization-P2} we note that
$(\mu_2)_*(f_2)=\mu_2\circ f_2$ apart of being contained in $W_2$ is
contained also in another component of $\Hom(\P^1,X)$ which dominates
$E_1$ and which contains morphisms from $W_1$ composed with a finite
morphism $\pi_d: \P^1\ra\P^1$ such that $f_1\circ\pi_d=\mu_2\circ f_2
$

We claim that $\Hom(\P^1,X)$ is not 4-connected at $\mu_2\circ f_2$,
more precisely that it is not connected after removing $\mu_2\circ
f_2\circ Aut(\P^1)$. Indeed, otherwise there would exist small
deformations of $\mu_2\circ f_2$ in $W_2$ which were contained in a
component of $\Hom(\P^1,X)$ different from $W_2$. Thus we would get a
component $W_3\subset \Hom(\P^1,X)$ dominating $E_2$. Lifting $W_3$ to
$\hat W_3\subset \Hom(\P^1,\P^2)$ we would get a component meeting $\hat
W_2$. This impossible because $\Hom(\P^1,\P^2)$ is smooth.

\end{pf}

\begin{prop}\label{normal-homeo}
The normalization map $\mu:\P^2\ra E$ is a homeomorphism in 
codimension $1$. 
\end{prop}

\begin{pf} Let $\Sigma=\{x\in E: |\mu^{-1}(x)|>1\}$ be the set over
which $\mu$ is not bijective. We argue by contradiction. Let us assume
that $\Sigma$ contains a 1-dimensional component $C_0$. Suppose that
$C_1,\dots, C_k$, with $k\geq 1$, are curves which are mapped via $\mu$
to $C_0$; by \ref{normalization-final} they are all rational curves. Let
$f_0:\P^1\ra C_0\subset X$ denote the normalization. Choose an
irreducible component $W_0\subset \Hom(\P^1,X)$ which contains $f_0$.
Then $W_0$ dominates $E$ and thus we can lift it to component
$W_1\subset \Hom(\P^1,\P^2)$ (recall that this means that $\mu_*(W_1)=
W_0$). If $f_1\in W_1$ is such $\mu_*(f_1)=f_0$ then $f_1$ maps $\P^1$
birationally onto a rational curve in $\P^2$, say $C_1$, and $\mu:
C_1\ra C_0$ is birational. Hence $k\geq 2$ and moreover $deg C_i =
deg(\mu_{|C_i})\cdot deg C_1$.

After re-numerating curves we may assume that $deg C_2$ is the smallest
among $degC_2,\dots, degC_k$. Let $W_2 \subset \Hom(\P^1,\P^2)$ be the set
parametrizing curves of degree $degC_2$ (recall that any component of
$\Hom(\P^1,\P^2)$ parametrizing morphisms of a given degree is smooth
and connected), so that the normalization $f_2: \P^1\ra C_2\subset \P^2$
is in $W_2$. Let $d=degC_2/degC_1$, then exists the unique degree $d$
cover $\pi_d:\P^1\ra\P^1$ such that $\mu\circ f_1\circ \pi_d= f_0\circ
\pi_d = \mu\circ f_2$ and clearly $f_1\circ\pi_d\ne f_2$.

Let $W_2'=\mu_*(W_2)$, then $W_2'$ is a connected component of
$\Hom(\P^1.X)$ (note that because of \ref{dominating-components} any
irreducible component of $\Hom(\P^1,X)$ dominates $E$ so can be lift up
to a component of $\Hom(\P^1,\P^2)$). We claim that $W_2'$ is not
4-connected (in the analytic topology !!) at $\mu\circ f_2$. Let us take
a small deformation of $f_2$, call it $h$, which is not contained in
$f_2\circ Aut(\P^1)$. Then the generic point of $h(\P^1)$ is in the set
where $\mu$ an isomorphism is and therefore $\mu\circ h$ lifts up only
to $h$, that is $\mu_*^{-1}(\mu_*(h))=h$. Thus, if we take a small
analytic neighbourhood $\ccU\ni \mu\circ f_2$, such that
$\mu^{-1}_*(\ccU)$ decomposes into disjoint small neighborhoods of the
inverse images $\mu^{-1}_*(\mu\circ f_2)\ni f_2,\ f_1\circ\pi_d$ then
$\ccU\setminus\mu\circ f_2\circ Aut(\P^1)$ will be disconnected.

\end{pf}


\subsection{Non-vanishing}

In this section we will use base-point-free techniques, due to
Kawamata, Koll\'ar. In \cite{AW-Duke} the method was developed for
studying local contractions, in the present situation we use
\cite{Mella} for reference as it is particularly applicable in our
case.

\begin{prop}\label{sect-O(1)}
There is a line bundle $\ccO_X(1)$ on $X$, such that $\nu^*\ccO_X(1)\cong
\ccO_{\P^2}(1)$.
\end{prop}

\begin{pf} The exponential sequences on $X$, $E$ and $\P^2$ together
with the vanishing of the higher cohomology groups of $\ccO_X$ and
$\ccO_{\P^2}$ give the following diagram coming from the normalization
$\mu: \P^2\ra E$ and the inclusion $i: E\ra X$.
\begin{eqnarray*}
\begin{array}{ccccc} \Pic X & \ra & \Pic E & \ra  & \Pic \P^2 \\ ||    
&     & \da    &      & ||        \\ \H^2(X,\Z) & \stackrel{i^*}{\ra} &
\H^2(E,\Z) & \stackrel{\mu^*}{\ra} & \H^2(\P^2,\Z)  \end{array}
\end{eqnarray*}
By \ref{picard=H2}, $i^*$ is an isomorphism. By
\ref{normal-homeo}, $\mu^*$ is an isomorphism. Therefore
$\Pic X \ra \Pic \P^2 $ is an isomorphism. \end{pf}

\begin{prop}\label{non-vanish-section}
After possibly shrinking $X$ the bundle $\ccO_X(1)$
has a section that doesn't vanish on $E$. 
If we denote its zero set by $X'$ then $X'$ is smooth.
\end{prop}

\begin{pf} As in \cite[Claim 3.1]{AW-Duke}, we find a $\Q$-divisor $D$ on $X$,
which is a pull back from $Y$ and such that $(X,D)$ is klt outside $E$
and lc at $E$ (see \cite[2.34]{KM} e.g. for definition).  Let
$W\subset X$ be a minimal center of the lc singularities of $(X,D)$
(see e.g. \cite{Mella} for a definition).  By construction $W\subset
E$. By Kawamata \cite[1.6]{Kawamata97}, see also
\cite[Thm.~1.12]{Mella}, $W$ is normal. Since $E$ may be assumed to be
non-normal, it follows that $\dim W\le 1$.  We can now apply
\cite[Lemma 2.2.ii]{Mella} (using Mella's notation we have in our case
$r=0$ and $\gamma=0$) to conclude that there exists a section of
$\ccO_X(1)$ which doesn't vanish identically on $W$. This finishes the
first part, to prove the second part it is enough to check smoothness
along $X'\cap E$. Let $p\in X'\cap E$ be a point. Since $E$ normalizes
to a $\P^2$ we find a line in $\P^2$ with image $C$ in $E$ such that
$p\in C$ and $C\not\subset X'$. By construction, the intersection
number $X'\cdot C$ in $X$ is $1$. This implies that both $X'$ and $C$
are smooth at $p$.
\end{pf}

Let $\f':X'\ra Y'$ be the normalized restriction of $\f$ to $X'$. Let
$\ccO_{X'}(1)$ be the restriction of $\ccO_X(1)$ to $X'$.  Note that by
adjunction, $K_{X'}=\ccO_{X'}(1)$. Let $F = E\cap X'\subset X$. In
$X'$ the set $F$ is $\f'$-exceptional. Note that $F$ is irreducible and it
is the image of a line in $\hat{E}=\P^2$.

\begin{prop}\label{two-cases}
In the above situation one of the two cases holds:\\ 
(i) $F$ is smooth.\\
(ii)  $\ccO_{X'}(1)$ has a section that does not vanish on $F$.\end{prop}

\begin{pf} 
Again, by \cite[Claim 3.1]{AW-Duke} we find a $\Q$-divisor $D$ on $X'$,
which is a pull back from $Y'$ and such that $(X',D)$ is lt outside
$F$ and lc inside $F$. Let $W$ be a minimal center of lc
singularities. Again, by Kawamata \cite{Kawamata97}, see \cite[Thm
1.12]{Mella}, $W$ is normal.  If $\dim W=1$, then it follows that $F$
is normal hence it is smooth and we are in case (i).

Let's suppose that $\dim W=0$.  We now wish to apply \cite[Lemma 2.3.iii]{Mella}. 
Using Mella's notation, we have $\gamma=0$, $r=-1$, $w=1$.
Therefore there exists a section of $\ccO_X'(1)$ that doesn't vanish
on $F$. This gives Case (ii). \end{pf}

Now we will discuss the cases.

\begin{prop}\label{case-one}
Suppose Proposition~\ref{two-cases}.i holds. 
Then $E$ is regular in codimension $1$.
\end{prop}

\begin{pf} Let $l\subset \P^2$ be the line such that $\mu(l)=F$. We
claim that $\mu: \P^2\ra E\subset X$ 
is an immersion along $l$. So we check the derivative
map $(T_{\P^2})_{|l}=T_l\oplus N_{l/\P^2}\ra 
(T_X)_{|F}=T_F\oplus N_{F/X}$. Since
$l\ra F$ is an isomorphism it is enough to verify $N_{l/\P^2}\ra
N_{F/X}$. Since however $F\subset X'$ thus we can consider the
composition $N_{l/\P^2} \ra N_{F/X}\ra (N_{X'/X})_{|F}$ which is an
isomorphism.
\end{pf}

The proof of the following proposition is postponed until the next section.

\begin{prop}\label{case-two}
Suppose Proposition~\ref{two-cases}.ii holds. 
Then $E$ is regular in codimension $1$.
\end{prop}

Assuming \ref{case-two}, in order to finish
the proof of Theorem~\ref{big-thm} it is enough to prove
the following result.

\begin{prop}\label{codim-1-regular-implies-big}
Suppose $E$ is regular in codimension $1$.
Let $\ccI\subset \ccO_{X}$ be the ideal sheaf 
defining $E$. Then\\
(i) $H^2(E,\ccI^r/\ccI^{r+1})=0$ for each $r\ge 1$, \\
(ii) $H^1(E,\ccO_E)=0,$ \\
(iii) $E$ is normal.
\end{prop}

\begin{pf}
 Let $\mu:\P^2\ra E$ be the normalization map. 

\noindent{\sl (i)} There is a natural surjective map 
$S^r (\ccI/\ccI^2)\ra \ccI^r/\ccI^{r+1}$. 
Therefore it is enough to show that $H^2(E,S^r (\ccI/\ccI^2))=0$. But 
since $E$ is Lagrangian, the sheaves $\mu_*S^rT_{E}$
and $S^r (\ccI/\ccI^2)$ are isomorphic over the smooth locus of $E$,
in particular they are isomorphic in codimension $1$. 
We note that any two sheaves that are isomorphic in codimension 
$1$ on have the same top cohomology class.
Therefore it is enough to show that
$\H^2(E,\mu_*S^r T_{\P^2})$ is zero. 
But this follows from the well know fact that 
$H^2(\P^2, S^rT_{\P^2})=0$ for all $r>0$.

\noindent{\sl (ii)} Let $E_r$ be the closed subscheme of $X$, 
defined by the ideal sheaf 
$\ccI^r$. The obstruction to lifting an element of 
$H^1(E_r,\ccO_{E_r})$ to $H^1(E_{r+1},\ccO_{E_{r+1}})$ lies
in $H^2(E,\ccI^r/\ccI^{r+1})$, which is zero. This implies
that the restriction map 
$\ilim_r  H^1(E_r,\ccO_{E_r}) \lra  H^1(E,\ccO_E)$
is surjective. By the Theorem on formal functions \cite[II.11.1]{Hartshorne}, 
the left hand side is 
isomorphic to $(R^1 f_* \ccO_{X})_0^{\wedge}$, 
which is zero by \ref{GRKKV-vanishing}.
Therefore $H^1(E,\ccO_E)=0.$

\noindent{\sl (iii)} We define the sheaf $\ccQ$ on $E$ to be the quotient
$$0\lra \ccO_E\lra \mu_* \ccO_{\P^2}\lra \ccQ\lra 0$$
Then $\ccQ$ has support precisely on the locus of points where $\mu$ is not an 
isomorphism, that is $\ccQ$ has only zero-dimensional 
support. Now take cohomology in above sequence to obtain 
$$0\lra \H^0(E,\ccO_E)\lra \H^0(E, \mu_*\ccO_{\P^2})\lra \H^0(E, \ccQ)\lra \H^1(E,\ccO_E).$$
We have $h^0(E,\ccO_E)=1$ and $h^0(E, \nu_*\ccO_{\P^2})=h^0(\P^2,\ccO_{\P^2})=1$.
By (ii), we have $h^1(E,\ccO_E)=0$. This implies
$h^0(E,\ccQ)=0$ and consequently $\ccQ=0$. 
This implies that $\P^2\cong E$.
\end{pf}


\subsection{A non-normal surface}

We assume from now on in this section that
Proposition~\ref{two-cases}.ii holds and moreover, in order to derive
the contradiction with \ref{case-two}, that $E$ is not regular in
codimension 1 so that in particular $\ccO_E(1)$ is not spanned.

\begin{lem}\label{sect-O(2)} In the situation of
\ref{two-cases}ii the line bundle $\ccO_X(2)$ is spanned. 
\end{lem}

\begin{pf} Since $\H^1(\ccO_X(1))=0$, the restriction map 
$\H^0(\ccO_X(2))\ra\H^0(\ccO_{X'}(2))$ is surjective. Therefore it is
enough to show that $\ccO_{X'}(2)$ is spanned. By  assumption,
$\ccO_{X'}(1)$ has a section that does not vanish  on $F$. Let
$X''\subset X'$ be its zero locus. Then $X''$ is  affine hence
$\ccO_{X''}(2)$ is spanned. Therefore it is enough to show that the
restriction map  $\H^0(\ccO_{X'}(2))\ra\H^0(\ccO_{X''}(2))$ is
surjective. For that it is enough to show that $\H^1(\ccO_{X'}(1))=0$,
which in turn follows from the exact sequence $\H^1(\ccO_X(1))\ra
\H^1(\ccO_{X'}(1))\ra\H^2(\ccO_X)$ and the vanishing of
$\H^1(\ccO_X(1))$ and $\H^2(\ccO_X)$. \end{pf}

Now we use the sections of $\ccO_X(1)$ and $\ccO_X(2)$ to analyse
the geometry of $E$. 

We pull-back sections of $\ccO_E(a)$ to sections of
$\ccO_{\P^2}(a)$. By \ref{sect-O(1)} and \ref{two-cases}ii,
$\ccO_E(1)$ has two sections, their pull-backs are denoted by $x,y\in
\H^0(\P^2,\ccO_{\P^2}(1))$. Therefore $\ccO_E(2)$ admits three
sections $x^2,xy,y^2\in \H^0(\P^2,\ccO_{\P^2}(2))$. Moreover, by
\ref{sect-O(2)}, $\ccO_E(2)$ admits a section $w\in
\H^0(\P^2,\ccO_{\P^2}(2))$, such that $(x^2,xy,y^2,w)$ do not vanish
simultaneously. Let $z\in \H^0(\P^2,\ccO_{\P^2}(1))$, such that
$[x,y,z]$ form a homogeneous coordinate system of $\P^2$. Then $w$ is
a quadratic polynomial in $x,y,z$. Let $\pi:\P^2\ra \P^3$ be the
morphism defined by $[x,y,z]\mapsto [x^2,xy,y^2,w]$; the image of
$\pi$ is the singular quadric cone $Q\subset\P^3$.  By the
construction the morphism $\pi$ factors through the normalization
$\mu$ so that we have $\eta:E\ra Q$ which is a $2:1$ covering.  Thus
$E$ is ``sandwiched'' $\P^2\ra E\ra Q$ between two known surfaces.

Let $e_0\in E$ be the unique base point of $\ccO_E(1)$.  Then $e_0$ is
a smooth point $E$ (the sections of $\ccO_E(1)$ pull-back to local
coordinates around the inverse image of $e_0$) and it is mapped to the
vertex $q_0$ of $Q$. Let us blow-up $Q$, respectively $E$ and $\P^2$,
along $q_0$, $e_0$ and $\mu^{-1}(e_0)$ to $Q'$, $E'$ and $\P'$ and
denote the exceptional curves by $A_Q$, $A_E$ and $A_\P$,
respectively.  The induced morphism of the blow-ups we denote by
$\mu'$, $\pi'$ and $\eta'$, respectively.  Then $Q'$ is rational ruled
with the exceptional curve $A_Q$ and the fiber of the ruling denoted
by $F_Q$. If $\ccO_{\P'}(1)$ denote the pullback of $\ccO_{\P^2}(1)$
then $\ccO_{\P^2}(1)=(\pi')^*(A_Q+2F_Q)$, moreover
$(\pi')^*(A_Q)=2A_\P$.  The branching divisor of $\pi':\P'\ra Q'$
consist of $A_Q$ and a curve $B_Q$, disjoint from $A_Q$, which is in
$|\ccO_{Q'}(A_Q+2F_Q)|$, the latter can be checked by adjunction.
Therefore $\pi'_*\ccO_{\P'}=\ccO_{Q'}\oplus\ccO_{Q'}(-A_Q-F_Q)$.  By
$B_\P\subset\P'$ let us denote the component of the ramification of
$\pi'$ which is over $B_Q$; it is a lift-up of a line $B\subset \P^2$.

The sequence $\P' \ra E'\ra Q'$ can be described in terms of
inclusions of $\ccO_{Q'}$ algebras:
$\ccO_{Q'}\ra\eta'_*\ccO_{E'}\ra\pi'_*\ccO_{\P'}$.  The trace splits
off the trivial factor $\ccO_{Q'}$ in $\pi'_*\ccO_{\P'}$ as well as in
$\eta'_*\ccO_{E'}$, so that we can write $\eta'_*\ccO_{E'}=
\ccO_{Q'}\oplus\ccJ$ for some rank 1 sheaf of $\ccO_{Q'}$ modules
which admits morphism $\ccJ\otimes\ccJ\ra \ccO_{Q'}$ coming from the
multiplication in $\ccO_{Q'}$ algebra structure.  Let $\ccS$ be the
reflexivisation of $\eta'_*\ccO_{E'}$. Since $\eta'_*\ccO_{E'}$ is
torsion free we have an inclusion $\eta'_*\ccO_{E'}\ra \ccS$ which
extends to
$$\ccO_{Q'}\ra\eta'_*\ccO_{E'}\ra\ccS\ra\pi'_*\ccO_{\P'}$$ 
Moreover,
we have the splitting $\ccS=\ccO_{Q'}\oplus\ccS_1$.  Being reflexive
of rank 1 the sheaf $\ccS_1$ is a line bundle and since it is the
reflexivisation of $\ccJ$ it admits the unique morphism
$\ccS_1\otimes\ccS_1\ra\ccO_{Q'}$ which gives a natural $\ccO_{Q'}$
algebra structure on $\ccS$. Thus we can set $S'=Spec_{Q'}\ccS$ and we
have a sequence of surjective morphism associated to the above
inclusions of sheaves
$$\P'\ra S'\ra E'\ra Q'$$

Since $\P'\ra E'$ is an isomorphism around $A_{\P}\ra A_E$ the image
of $A_\P$ in $S'$ is a $(-1)$-curve $A_S$ 
which can be blow-down to a smooth
point on a surface $\beta: S'\ra S$. Thus we get a sequence
$$\P^2\ra S\ra E\ra Q$$
and $S$ is a ``partial normalization''.
By the construction the resulting morphism $\alpha: 
S\ra E$ is an isomorphism in codimension 1
and $S$ satisfies the Serre's condition $S_2$
hence it is the $S_2$-ization of $E$
c.f.~\cite{Reid}.
By $\ccO_S(1)$
let us denote the pull-back of $\ccO_E(1)$ and by $B_S$ the (reduced)
image of $B$, the ramification divisor of $\pi$.
 
\begin{lem}\label{S2-isation}
The surface $S$ is a del Pezzo non-normal surface: it has
locally complete intersection singularities (hence it is Cohen-Macaulay),
its non-normal locus is $B_S$  and its dualising sheaf is $\ccO_S(-1)$.
\end{lem}

\begin{pf}
The surface $S'$ is produced  as divisor in the total space of the 
line bundle
$\ccS_1$ given by the section coming from the multiplication
$\ccS_1\otimes\ccS_1\ra\ccO_{Q'}$ hence $S$
has only locally complete intersection singularities.
Since outside $B$ and $A_S$ the surface $S'$ was obtained by
factorizing 
a local analytic isomorphism coming from $\P'\ra Q'$
it is smooth; the smoothness of $S$ at the image of $A_S$ comes from the
construction. Thus it remains to compute the dualising sheaf of
$S$. For this we compute the line bundle $\ccS_1$:
by the inclusion $\ccS_1\ra\ccO_{Q'}(-A_Q-F_Q)$ we know that
$\ccS_1=\ccO(-A_Q-F_Q-D)$ where $D$ is an effective and non-zero divisor
(we assume that $E$ is not regular in codimension 1 !)
such that $D\cap A_Q=\emptyset$ and thus $D\in |d(A_Q+2F_Q)|$
where $d>0$.
On the other hand by the vanishing \ref{top-vanishing}
we know that $H^2(S,\ccO_S)=0$ hence also $H^2(Q',\ccS_1)=0$.
Comparing this with the previous observation and using duality
we get $H^0(Q',\ccS_1^*\otimes K_{Q'})=H^0(Q',\ccO_{Q'}((d-1)A_Q+(2d-3)F_Q)=0$
Therefore $d=1$. Now, denoting by $\gamma$ the morphism $S'\ra Q'$
we can compute $K_{S'}$ by adjunction:
$$K_{S'}=\gamma^*(K_{Q'}+2A_Q+3F_Q)=\gamma^*(-F_Q)= 1/2\cdot\gamma^*A_Q - 
1/2\cdot\gamma^*(A_Q+2F_Q)=A_S-\beta^*(\ccO_S(1))$$ 
Hence $K_S=\ccO_S(-1)$.
\end{pf}

The surface $S$ can be found in the list of \cite{Reid}, in particular
it can be written explicitly as the hypersurface $(z^2=y^3)$
in the weighted projective space $\P(1,1,2,3)$ with homogeneous
coordinates $[x_1,x_2,y,z]$. This allows to compute
its sheaf of differentials.

\begin{lem}\label{diff-S}
Let $S$ be the above surface with the non-normal locus
at $B_S\iso\P^1$. Then $\Omega_{S}|_{B_S}=\ccO_{\P^1}(-2)
\oplus\ccO_{\P^1}(-2 )\oplus\ccO_{\P^1}(-3)$
\end{lem}

\begin{pf}
In the above situation, let $\P$ denote the weighted projective space 
$\P(1,1,2,3)$ The line $B_S\subset\P$
is then given by the equations $y=z=0$ and it is contained in the smooth
locus of $\P$. Thus we can verify that 
$\Omega_{\P}|_{B_S}=\ccO_{\P^1}(-2)
\oplus\ccO_{\P^1}(-2 )\oplus\ccO_{\P^1}(-3)$.
On the other hand restricting the exact sequence
$I_{S/\P}/I_{S/\P}^2\ra \Omega_{\P}|_S\ra \Omega_S\ra 0$
to the line $B_S$ we get the isomorphism
$\Omega_{\P}|_{B_S}\ra \Omega_S|_{B_S}$.
\end{pf}

\noindent{\it Proof of \ref{case-two}}
We argue by contradiction. If $E$ is not normal in codimension
1 then the partial normalization
$\alpha: S\ra E$ is an isomorphism 
in codimension $1$, 
and the induced morphism of differentials
$\alpha^*\Omega_E|_{B_S}\ra \Omega_{S}|_{B_S}$ is an 
isomorphism at the generic point of the line $B_S$.
Moreover we have a 
surjective map $\Omega_X\ra \Omega_E$ coming from the inclusion
$E\ra X$.
Therefore we obtain a map $\alpha^*\Omega_X\ra \Omega_{E'}|L$ 
which is generically surjective. 
By the previous lemma it follows that $\alpha^*\Omega_X\cong
\sum_{i=1}^4\ccO(a_i)$ with $a_1\le a_2\le a_3\le a_4$ and 
such that $a_3<0$. But since $X$ is symplectic, we have 
$a_1=-a_4$ and $a_2=-a_3$. Contradiction.

\medskip

This finishes the proof of Theorem~\ref{big-thm}.


\bigskip
\hrule
\bigskip
\noindent
Jan Wierzba: Berlin, Germany
\par\noindent{\tt Jan.Wierzba@ATKearney.com}
\par\noindent
Jaros{\l}aw A. Wi\'sniewski: Institute of Math., Warsaw University, Banacha 2, Warsaw, Poland
\par\noindent{\tt
jarekw@mimuw.edu.pl}

\end{document}